\documentclass[12pt]{amsart}
\usepackage{pdfsync}
\usepackage{amsmath,amsthm,amsfonts,amssymb}
\usepackage[all]{xy}
\usepackage{amssymb, amsthm, amsmath}
\usepackage[english]{babel}
\usepackage{amsfonts}
\numberwithin{equation}{section}

\DeclareMathOperator{\id}{id}

\newcommand{\ra}{\rightarrow}

\newcommand{\C}{\mathbb C}\newcommand{\mS}{\mtc{S}}

\newcommand{\ot}{\otimes}
\newcommand{\mtc}{\mathcal}

\newcommand{\Lam}{\Lambda}
\newcommand{\Irr}{\mtr{Irr}}

\newcommand{\al}{\alpha}
\newcommand{\eps}{\epsilon}
\newcommand{\nhs}{normal Hopf subalgebra}
\newcommand{\sha}{semisimple Hopf algebra}

\newtheorem{lemma}[equation]{Lemma}
\newtheorem{thm}[equation]{Theorem}
\newtheorem{prop}[equation]{Proposition}
\newtheorem{defn}[equation]{Definition}
\newtheorem{cor}[equation]{Corollary}
\newtheorem{rem}[equation]{Remark}

\newcommand{\uw}{\uparrow}

\newcommand{\ch}{\chi}
\newcommand{\mtr}{\mathrm}
\newcommand{\kr}{\mathrm{ker}\;}
\newcommand{\bn}{\begin}

\newcommand{\ncm}{\newcommand}
\ncm{\np}{\newpage}
\ncm{\ebl}{\end{thebibliography}}
\ncm{\bbl}{\begin{thebibliography}}
\ncm{\chd}{_{ _{\ch}}}
\ncm{\ald}{_{ _{\al}}}
\ncm{\cP}{\mathcal{P}}
\ncm{\ei}{e_i}
\ncm{\eij}{e_{i,\;j}}
\ncm{\bt}{\begin{thm}}

\ncm{\bdef}{\begin{defn}
}

\ncm{\edf}{\end{defn}}
\ncm{\et}{\end{thm}}
\ncm{\bc}{\begin{cor}}
\ncm{\bl}{\begin{lemma}}
\ncm{\el}{\end{lemma}}
\ncm{\bpf}{\begin{proof}}
\ncm{\epf}{\end{proof}}
\ncm{\ec}{\end{cor}}
\ncm{\er}{\end{rem}}
\ncm{\br}{\begin{rem}}
\ncm{\bp}{\begin{prop}}
\ncm{\ep}{\end{prop}}
\ncm{\bd}{\begin{document}}
\ncm{\ed}{\end{document}}
\ncm{\beq}{\begin{equation}}
\ncm{\beqn}{\begin{equation*}}
\ncm{\eeq}{\end{equation}}
\ncm{\eeqn}{\end{equation*}}
\ncm{\bea}{\begin{eqnarray}}
\ncm{\eea}{\end{eqnarray}}
\ncm{\beanon}{\begin{eqnarray*}}
\ncm{\eeanon}{\end{eqnarray*}}
    \usepackage{fancyhdr}                    

\title[Semisimple Hopf algebras]
{Categorical Hopf kernels  and representations of semisimple Hopf algebras}
\author{Sebastian  Burciu}
\address{Inst.\ of Math.\ ``Simion Stoilow" of the Romanian Academy
P.O. Box 1-764, RO-014700, Bucharest, Romania }
\email{sebastian.burciu@imar.ro}
\date{\today}
\begin{document}
\maketitle
\begin{abstract}
In the category of semisimple Hopf algebras  the Hopf kernels introduced by Andruskiewitsch and Devoto in \cite{AD} coincide with kernels of representation as introduced in \cite{Bker}. Some new results concerning the normality of kernels are also presented. It is proven that the property for Hopf algebras to have all kernels normal Hopf subalgebras is a selfdual property.
\end{abstract}
\section{Introduction}
Semisimple Hopf algebras were intensively studied in the last twenty years. Many properties from finite groups are extended to the more general setting of semisimple Hopf algebras.

In order to work inside the category of Hopf algebras in \cite{AD} the authors introduced the notion of Hopf kernel of a morphism between Hopf algebras. Using the exponent properties of semisimple Hopf algebras in \cite{Bker} a notion of kernel of a representation of a semisimple Hopf algebra was proposed. This notion extends the notion of kernel of finite group representations. In this paper we prove that the two notions of kernels coincide. It will be proven in Proposition \ref{newdescr} that the Hopf kernel of a morphism between semisimple Hopf algebras is the kernel of a certain representation. Theorem \ref{coinc} shows that the converse of this fact is also true, any kernel of a representation is also a Hopf kernel in the sense defined in \cite{AD}.

Although many properties of kernels were transferred from groups to this more general setting of semisimple Hopf algebras there are still some unanswered questions in this direction. For example the normality of kernels of representations was proven in \cite{Bker} with an additional assumption, that of centrality of the character in the dual Hopf algebra. It is not known yet if this additional assumption is necessary. In this paper we study Hopf algebras where all these kernels are normal Hopf subalgebras. We say that a Hopf algebra has property $(N)$ if any kernel $H_{\ch}$ is a normal Hopf subalgebra of $H$ for any $ \ch \in \Irr(H)$. It will be shown in Theorem \ref{sfdual} that property $(N)$ is a self dual notion, $H$ has property $(N)$ if and only if $H^*$ has property $(N)$. We also give new information on these kernels in terms of the central characters described in \cite{Zh}.

We work over the base field $\C$ and all Hopf algebra notations are those from \cite{Montg}. We drop the sigma symbol in Sweedler's notations for comultiplication.
\section{The kernel of a representation as a Hopf kernel}\label{kers}

Through all this paper $H$ is a semisimple Hopf algebra over $\C$. It follows that $H$ is also cosemisimple and $S^2=\id_H$. The set of irreducible characters of $H$ is denoted by $\Irr(H)$ and this is a base of the character algebra $C(H)\subset H^*$. Moreover $C(H)$ is a semisimple algebra \cite{Z}.

To any irreducible character $d \in \mtr{Irr}(H^*)$ is associated a simple comatrix coalgebra $C=\mathbb{C}<x_{ij}>_{1\leq i,j \leq q}$ as in \cite{Lar}.

If $\Lam$ is the idempotent integral of $H$ it follows that $ \dim_{\C}(H)\Lam$ is the regular character of $H^*$, that is:
\begin{equation}\label{Lam}
    \dim_{\C}(H)\Lam=\sum_{d \in \Irr(H^*)}\eps(d)d.
\end{equation}

Let $H$ be a semisimple Hopf algebra over $\mathbb{C}$ and $M$ be an $H$-module affording the character $\ch$. Proposition 1.2 from \cite{Bker} shows that $|\ch(d)|\leq \eps(d)\ch(1)$. In fact the equivalence of the following assertions follows from the same Proposition.

\bn{prop}\label{equival}( see \cite{Bker} Remark 1.3)
Let $H$ be a semisimple Hopf algebra over $\mathbb{C}$ and $M$ be an $H$-module affording the character $\ch$.
\begin{enumerate}
    \item $\ch(d)=\eps(d)\ch(1)$.
    \item $\ch(x_{ij})=\delta_{ij}\ch(1)$ for all $i,j$.
    \item $dm=\eps(d)m$ for all $m\in M$.
    \item $x_{ij}m=\delta_{ij}m$ for all $i,j$ and $m \in M$.
\end{enumerate}
\end{prop}

The kernel $H_{ _{M}}$(or $H_{ _{\ch}}$) is defined as follows (see \cite{Bker}). Let  $\mtr{ker}_H(\ch)$ be the set of all irreducible characters $d \in \mtr{Irr}(H^*)$ which satisfy the equivalent conditions above. It can be proven that this set is closed under multiplication and $``^*"$ and therefore it generates a Hopf subalgebra $H_{ _M}$ (or $H_{ _{\ch}}$) of $H$ which is called the kernel of the representation $M$ \cite{NR}.

\begin{rem}\label{st}\end{rem}
\begin{enumerate}

\item \label{prod}
For later use let us notice that $\mtr{ker}_H(\chi) \subset \mtr{ker}_H(\chi^n)$ for all $n \geq 0$. This is item 1 of Remark 1.5 from \cite{Bker}.  It follows that $\cap_{n \geq 0}\mtr{ker}_H(\chi^n)=\mtr{ker}_H(\chi)$ which can also be written as $\cap_{n \geq 0}H_{ _{M^{\ot n}}}=H_{ _M}$.

\item \label{conj}
We also need the following result proven in \cite{NR'}. Suppose that $d$ is a character of $H^*$ and $\ch$ a character of $H$. Then $\chi(d^*)=\overline{\chi(d)}$.

\item \label{equivz} Suppose that $\ch_i \in \Irr(H)$ and $\chi=\sum_{i=1}^{s}m_i\chi_i$ where $m_i>0$. Then it is easy to see that
\begin{equation*}
\mtr{ker}_H(\ch)=\cap_{i=1}^r\mtr{ker}_H(\chi_i)
\end{equation*}
\end{enumerate}

\bn{prop}\label{larg}
Let $H$ be a semisimple Hopf algebra over $\mathbb{C}$ and $M$ be a representation of $H$ with character $\ch$. Consider the subalgebra of $H$ given by $$\mS_{ _M}=\{h \in H\;|\; hm=\eps(h)m\; \text{for all} \;\;m \in M\}.$$
Then $H_{ _M}$ is the largest Hopf subalgebra of $H$ contained in $\mS_{ _M}$.
\end{prop}

\bn{proof}
By the definition of the kernel it is clear that $H_{ _M} \subset \mS_{ _M}$. Item \ref{conj} of Remark \ref{st} implies that if $C \subset \mS_{ _M}$ then $S(C) \subset \mS_{ _M}$. Previous Porposition also shows that if a subcoalgebra $C$ is included in $\mS_{ _M}$ then $C$ is also included in $H_{ _M}$. It is easy to see that if $C$ and $D$ are two subcoalgebras included in $\mS_{ _M}$ then the product coalgebra $CD$ is also included in $H_{ _M}$. It follows that the largest Hopf subalgebra included in $\mS_{ _M}$ it is the sum of all subcoalgebras included in $\mS_{ _M}$ and therefore this Hopf subalgebra is also included in $H_{ _M}$. Thus this Hopf subalgebra coincides with $H_{ _M}$
\end{proof}

\bn{cor}
Let $H$ be a semisimple Hopf algebra over $\mathbb{C}$ and $M$ be a representation of $H$. Then
$H_{ _M}$ is the largest Hopf subalgebra $K$ of $H$ such that $ K^+ \subset \mtr{Ann}_H(M)$.
\end{cor}

\bn{proof}
It is easy to see that for any Hopf subalgebra $K$ of $H$ one has $K^+ \subset \mtr{Ann}_H(M)$ if and only if $K \subset H_{ _M}$
\end{proof}
\bn{rem}\end{rem}\label{2star}
Let $K$ be a normal Hopf subalgebra of $H$ and $L:=H//K$ be the quotient Hopf subalgebra. From previous Corollary it follows that $\Irr(L)=\{\ch \in \Irr(H)\;|\; H_{\ch} \supset K\}$.
\subsection{Hopf kernel of a Hopf algebra map}
In this subsection it will be shown that any Hopf kernel is a kernel of representation. The converse it will be proven in Theorem \ref{coinc}.

Recall the Hopf kernel of a Hopf map defined in \cite{AD}. If $f :A \ra B$ is a Hopf algebra map then the Hopf kernel of $f$ is defined as follows:

\begin{equation}\label{Hopfk}
    HKer(f)=\{a\in H\;|\; a_1\ot \pi(a_2) \ot a_3 =a_1\ot \pi(1) \ot a_2\}.
\end{equation}

It was proven in \cite{AD} that the Hopf kernel is a Hopf subalgebra of $H$.
\bn{prop}\label{newdescr}
Let $I$ be a Hopf ideal of $H$ and $\pi : H \ra H/I$ be the canonical Hopf projection. Then

\begin{equation*}
    \mS_{ _{H/I}}=\{h \in H \;\; | \pi(h) =\eps(h) 1\}.
\end{equation*}

Regarding $H/I$ as $H$-module it follows that $\mtr{HKer}(\pi)=H_{ _{H/I}}$.
\end{prop}

\bn{proof}
It is straightforward to verify the formula for $\mS_{ _{H/I}}$. It is also clear that $\mtr{HKer}(\pi) \subset \mS_{ _{H/I}}$.  Since $\mtr{HKer}(\pi)$ is a Hopf subalgebra of $H$ Proposition \ref{larg} implies that $\mtr{HKer}(\pi)\subset H_{ _{H/I}}$. On the other hand it easy to check using the equivalencies from Proposition \ref{equival} that $H_{ _{H/I}}\subset \mtr{HKer}(\pi)$. Therefore the equality $\mtr{HKer}(\pi)=H_{ _{H/I}}$ holds.
\end{proof}

\subsection{Description of the categorical Hopf kernel as the kernel of a representation}
Let $M$ be an $H$-module with character $\ch$. Consider $$I_M=\cap_{n=0}\mtr{Ann}_H(M^{\ot n}).$$ Then $I_M$ is a Hopf ideal and it is the largest Hopf ideal contained in the annihilator $\mtr{Ann}_H(M)$ \cite{KSZ}. Therefore $B=H/I_M$ is a Hopf algebra and one has a canonical projection of Hopf algebras $$\pi:H \ra B.$$

\bn{prop}Let $H$ be a semisimple Hopf algebra over $\mathbb{C}$ and $M$ be an $H$-module. Using the above notations it follows that  $H_{ _M}=H_{ _B}$ where $B$ is regarded as $H$-module via $\pi$.
\end{prop}

\bn{proof}
The representations of the Hopf algebra $H/I_M$ are exactly the representations that are constituents in some power $M^{\ot n}$ of $M$. Thus one has $H_{ _B}= \cap_{n \geq 0}H_{ _{M^{\ot n}}}$ by item \ref{equivz} of Remark \ref{st}. Also item \ref{prod} of the same Remark \ref{st} implies that $H_{ _M}=H_{ _B}$.
\end{proof}

The next theorem gives the characterization of the kernel of a representation as a Hopf kernel.

\bn{thm}\label{coinc}
Let $H$ be a semisimple Hopf algebra and $M$ be a representation of $H$. Let as above $\pi:H \ra H/I_M$ be the canonical projection. Then
\begin{equation*}
    H_{ _M}=\mtr{HKer}(\pi).
\end{equation*}
\end{thm}

\bn{proof}
By previous Proposition it is enough to show that $$\mtr{HKer}(\pi)=H_{ _B}. $$

It is easy to see that $H_{ _B} \subset  \mtr{HKer}(\pi)$. Indeed, if $h\in H_B$ then $\sum h_1 \ot \pi(h_2) \ot h_3=\sum h_1\ot \eps(h_2)1 \ot h_3=\sum h_1\ot 1 \ot h_2$ since $\pi(h)=\eps(h)1$ for all $h \in H_{ _B} $ .

One can also see that $\mtr{HKer}(\pi) \subset H_{_B}$. Indeed if $h \in \mtr{HKer}(\pi)$ then $\sum h_1 \ot \pi(h_2) \ot h_3=\sum h_1\ot 1\ot  h_2.$ Applying $\eps \ot id \ot \eps$ to this identity it follows that $\pi(h)=\eps(h)1$ and Lemma \ref{newdescr} implies $\mtr{HKer}(\pi) \subset H_{_B}$. Since $\mtr{HKer}(\pi)$ is a Hopf subalgebra it follows from Proposition \ref{larg} that  $\mtr{HKer}(\pi) \subset H_{_B}$. Thus $\mtr{HKer}(\pi)=H_{ _B}. $
\end{proof}

\section{Hopf algebras with all kernels normal}

In this section we describe some properties of Hopf subalgebras with all kernels $H_{\ch}$ normal Hopf subalgebras.
\subsection{Central characters in the dual Hopf algebra}\label{centralch}Let $H$ be finite dimensional semisimple Hopf algebra over $\mathbb{C}$. Define a central subalgebra of $H$ by $\hat{\mtr{Z}}(H):=\mtr{Z}(H) \bigcap C(H^*)$. It is the subalgebra of $H^*$-characters which are central in $H$. Let $\hat{\mtr{Z}}(H^*):=\mtr{Z}(H^*) \bigcap C(H)$ be the dual concept, the subalgebra of $H$-characters which are central in $H^*$.

Let $\phi:H^* \ra H$ given by $\;f \mapsto f\rightharpoondown\Lambda_{ _H} $ where $f\rightharpoondown\Lambda_{ _H}=f(S({\Lambda_{ _H}}_{ _1})){\Lambda_{ _H}}_{ _2}$. Then $\phi$ is an isomorphism of vector spaces \cite{Montg}.

\bn{rem} \label{formulae}
It can be checked that $\phi(\xi_d)=\frac{\eps(d)}{|H|}d^*$ and $\phi^{-1}(\xi_{ _\ch})=\ch(1)\ch$ for all $d \in \mtr{Irr}(H^*)$ and $\ch \in \mtr{Irr}(H)$ (see for example \cite{Montg}). Here $\xi_{\ch} \in H$ is the central primitive idempotent of $H$ associated to $\ch$. Dually, $\xi_d \in H^*$ is the central primitive idempotent of $H^*$ associated to $d \in \Irr(H^*)$.
\end{rem}

The following description of $\hat{\mtr{Z}}(H^*)$ and $\hat{\mtr{Z}}(H)$ was given in \cite{Zh}. Since $\phi(C(H))=\mtr{Z}(H)$ and $\phi(\mtr{Z}(H^*))=C(H^*)$ it follows that the restriction $$\phi| _{ _{\hat{\mtr{Z}}(H^*)}}:\hat{\mtr{Z}}(H^*)\ra\hat{\mtr{Z}}(H)$$ is an isomorphism of vector spaces.

Since $\hat{\mtr{Z}}(H^*)$ is a commutative semisimple algebra it has a vector space basis given by its primitive idempotents. Since $\hat{\mtr{Z}}(H^*)$ is a subalgebra of $\mtr{Z}(H^*)$ each primitive idempotent of $\hat{\mtr{Z}}(H^*)$ is a sum of primitive idempotents of $\mtr{Z}(H^*)$. But the primitive idempotents of $\mtr{Z}(H^*)$ are of the form $\xi_d$ where $d \in \mtr{Irr}(H^*)$. Thus, there is a partition $\{\mathcal{Y}_j\}_{j\in J}$ of the set of irreducible
characters of $H^*$ such that the elements $(e_j)_{j \in J}$ given by $$e_j=\sum_{d \in \mtc{Y}_j}\xi_d$$ form
a basis for $\hat{\mtr{Z}}(H^*)$. Since $\phi(\hat{\mtr{Z}}(H^*))=\hat{\mtr{Z}}(H)$ it follows that $\widehat{e}_j:=|H|\phi(e_j)$ is a basis for $\hat{\mtr{Z}}(H)$. Using the first formula from Remark \ref{formulae} one has

\begin{equation}\label{decp}
    \widehat{e}_j=\sum_{d \in \mtc{Y}_j}\eps(d)d^*.
\end{equation}

Proposition 3.3 of \cite{Bker} shows that kernels of central characters are normal Hopf subalgebras. Thus with the above notations $\kr_{H^*}(\widehat{e}_j)$ is a normal Hopf subalgebra of $H^*$.
\bn{rem}\label{st2}\label{part} By duality, the set of irreducible
characters of $H$ can be partitioned into a finite collection of subsets
$\{\mathcal{X}_i\}_{i\in I}$ such that the elements $(f_i)_{i \in I}$ given by

\begin{equation}\label{3star}
f_i=\sum_{\ch \in \mtc{X}_i}\ch(1)\ch
 \end{equation}

form a $\mathbb{C}$-basis for $\hat{\mtr{Z}}(H^*)$. Then the elements $\phi(f_i)=\sum_{\ch \in
\mtc{X}_i}\xi_{ _\ch}$ are the central orthogonal primitive idempotents of $\hat{\mtr{Z}}(H)$ and therefore they form a basis for this space. Clearly $|I|=|J|$.
\end{rem}
Write $d \sim d'$ if both appear in the same minimal central character $\widehat{e}_j$ of $H$. Clearly $\sim$ is an equivalence relation with equivalence classes $\mtc{Y}_j$.

For an irreducible character $d \in \mtr{Irr}(H^*)$ let $N(d)$ be the smallest Hopf subalgebra of $H$ containing $d$. This always exists since intersection of normal Hopf subalgebras is always a \nhs.

\bn{prop}\label{eqd} Suppose that $d ,d' \in \Irr(H^*)$ with $d \sim d'$.
Then $N(d)=N(d')$.
\end{prop}

\bn{proof}
Since $N(d)$ is central it follows from \cite{masnr} that the idempotent integral $\Lam$ of $N(d)$ is central in $H$. Since $N(d)$ is a semisimple Hopf algebra Equation \ref{Lam} implies that $\Lam$ is a scalar multiple of the sum $\sum_{e \in \Irr(N(d)^*)}\eps(e)e$. But $\Irr(N(d)^*)\subset\Irr(H^*)$ and decomposition \ref{decp} of central characters of $H^*$ shows that $d' \in N(d)$. Therefore $N(d') \subset N(d)$. Symmetry implies that $N(d')=N(d)$.
\end{proof}
\bn{rem}\label{bij}(see Remark 2.3 of \cite{Bker}.)

1) Suppose that $K$ is a normal Hopf subalgebra of $H$ and let $L=H//K$ be the quotient Hopf algebra of $H$ via $\pi:H\ra L$. Then $\pi^* :L^* \ra H^*$ is an injective Hopf algebra map. It follows that $\pi^*(L^*)$ is a normal Hopf subalgebra of $H^*$. Moreover ${(H^*//L^*)}^*\cong K$ as Hopf algebras.

2) There is a bijection between normal Hopf subalgebras of $H$ and $H^*$. To any normal Hopf subalgebra $K$ of $H$ one associates $(H//K)^*$ as \nhs \; of $H^*$. Conversely to any $L$, a \nhs \;in $H^*$, one associates $(H^*//L)^*$ as normal Hopf subalgebra of $H$. The fact that these maps are inverse one to the other follows from the previous item of this Remark.

\end{rem}

\subsection{Hopf algebras with all kernels normal}
\bn{defn} We say that a semisimple Hopf algebra $H$ has property $(N)$ if and only if $H_{\ch}$ is a normal Hopf subalgebra of $H$ for all irreducible characters $\ch \in \Irr(H)$.
\end{defn}
Dually $H^*$ has property $(N)$ if and only if $H^*_d$ is normal for any irreducible character $d \in \Irr(H^*)$.

For a semisimple Hopf algebra $H$ denote by $t_H \in H^*$ the idempotent integral of $H^*$ It follows as in Equation \ref{Lam} that $\dim_{\C}(H)t_H$ is the regular character of $H$. The following theorem from \cite{Bker} will be used in the sequel.

\begin{thm}\label{main}
Let $H$ be a finite dimensional semisimple Hopf algebra. Any normal Hopf subalgebra $K$ of $H$ is the kernel of a character which is central
in $H^*$. More precisely, with the above notations one has:
\bn{equation*}
K=H_{{|L|\pi^*(t_{L})}}
\end{equation*}
where $L=H//K$ and $t_L$ is the idempotent integral of $L$.
\end{thm}

From the proof of this theorem also follows that the regular character of $L$ is $\dim_{\C}(L)t_L=\eps_K\uw^H_K$.

\bn{thm}
Let $H$ be a \sha.
If $H^*$ has property $(N)$ then $\mtr{ker}\;d=\mtr{ker}\;d'$ for all $d \sim d'$.
\end{thm}

\bn{proof}
Suppose $H^*_d$ is a normal Hopf subalgebra of $H$. Then by Remark \ref{bij} there is some normal Hopf subalgebra $K$ of $H$ such that $H^*_d=(H//K)^*$. From the definition of the kernel  and Remark \ref{2star} one has that $$\mtr{ker}_{H^*}(d)=\{\ch\;| \mtr{ker}\;\ch \supset K\}.$$

We claim that $d \in K$. Indeed as above the regular character of $H//K$ is $\eps_K\uw_K^H$. Note that $d \in \mtr{ker}_H(\ch)$ if and only if $\ch \in \mtr{ker}_{H^*}(d)$. This implies that $d \in H_{\eps_K\uw_K^H}$. But Corollary 2.5 of \cite{Bker} implies that $H_{\eps_K\uw_K^H}=K$.

Since $N(d)$ is the minimal normal Hopf subalgebra containing $d$ one has that $N(d) \subset K$. On the other hand since $\Lam_K$ is central in $H$  decomposition \ref{decp} of central characters of $H^*$ implies $d' \in K$. Therefore $\mtr{ker}\;d' \supset \mtr{ker}\;d$. By symmetry one obtains $\mtr{ker}\;d' \subset \mtr{ker}\;d$ and thus the equality $\mtr{ker}\;d' = \mtr{ker}\;d$.
\end{proof}

\bn{cor}
Let $H$ be a \sha . Then $H^*$ has property $(N)$ if and only if $\mtr{ker}\;d=\mtr{ker}\;d'$ for all $d \sim d'$.
\end{cor}

\bn{proof}
We have already shown that if $H^*$ has $(N)$ then $\mtr{ker}\;d=\mtr{ker}\;d'$ for all $d \sim d'$. The converse follows by item \ref{equivz} of Remark \ref{st} since $\mtr{ker}_{H^*}(d)= \cap_{d ' \in \mtc{Y}_j}\mtr{ker}_{H^*}(d')= \kr_{H^*}(\widehat{e}_j)$ which is normal by Proposition 3.3 of \cite{Bker}.
\end{proof}

Now we can prove the main result of this subsection.
\bn{thm}\label{sfdual}
Let $H$ be a \sha . Then $H^*$ has property $(N)$ if and only if $H$ has property $(N)$.
\end{thm}

\bn{proof}
We show that if $H^*$ has $(N)$ then $H$ has also $(N)$. Write $\ch \sim \ch'$ if both characters appear in a minimal decomposition of a central characters $f_i$ of $H$ from Remark \ref{part}. Suppose $\ch \sim \ch'$ and let $d \in \mtr{ker}_H(\ch)$. Since $H^*$ has $(N)$ as before we have that $H^*_d=(H//K)^*$ for some $K$, a \nhs\; of $H$. Then Theorem \ref{main} implies that $\ch$ is a constituent of $\eps_K\uw_K^H$. Since by the same Theorem $\eps_K\uw_K^H$ is central in $H^*$ it follows from decomposition \ref{3star} of central characters of $H$ that $\ch'$ is also a constituent of $\eps_K\uw_K^H$. Therefore also $d \in \mtr{ker}\;\ch'$. This shows $ \mtr{ker}\;\ch\subset \mtr{ker}\;\ch'$. By symmetry one has the equality $ \mtr{ker}\;\ch= \mtr{ker}\;\ch'$. Previous Corollary shows that $H$ has also property $(N)$.
\end{proof}

\bn{prop}
Suppose that $H$ is a \sha\; with property $(N)$ and let $d \in \Irr(H^*)$. In this situation $H^*_d=(H//N(d))^*$.
\end{prop}

\bn{proof}
It is enough to show that $\mtr{ker}_{H^*}(d)=\Irr(H//N(d))$. Clearly $\Irr(H//N(d))\subset\mtr{ker}_{H^*}(d) $ by Remark \ref{2star}. Suppose now that $\ch \in \mtr{ker}_{H^*}(d)$. Since $H$ has $(N)$ it follows that $H_{\ch}$ is a normal Hopf subalgebra of $H$. Since $d \in \mtr{ker}_H(\ch)$ definition of $N(d)$ shows that $N(d)\subset H_{\ch}$. Thus $\ch \in \Irr(H//N(d))$.
\end{proof}
\subsection*{Acknowledgments}
This work was supported by the strategic grant POSDRU/89/1.5/S/58852, Project "Postdoctoral programme for training scientific researchers" cofinanced by the European Social Found within the Sectorial Operational Program Human Resources Development 2007 - 2013.
\bibliographystyle{amsplain}
\bibliography{categkernels}
\end{document}